\documentclass[12pt]{amsart}
\usepackage{amssymb}
\usepackage{times,a4wide}

\newtheorem{theorem}{Theorem}[section]
\newtheorem{definition}[theorem]{Definition}
\newtheorem{proposition}[theorem]{Proposition}
\newtheorem{corollary}[theorem]{Corollary}
\newtheorem{lemma}[theorem]{Lemma}

\def \proof {\noindent {\bf Proof.}\ \ }
\def \qed {{\mbox{}\nolinebreak\hfill\rule{2mm}{2mm}\par\medbreak} }
\def \remark {\noindent {\bf Remark.}\ \ }

\def \a {\alpha}

\def \e {\varepsilon}
\def \d {\delta}

\def \l {\lambda}

\def \la {\langle}
\def \ra {\rangle}
\def \o {\omega}
\def \s {\sigma}

\def \w {\omega}

\def \R {\mathbb{R}}

\def \rank {{\rm rank}}
\def \tr {{\rm tr}}

\def \Z {\mathcal{Z}}
\def \C {\mathbb{C}}
\def \N {\mathbb{N}}

\newcommand{\He}{\mathcal{H}}
\def \S {\mathbb{S}}
\newcommand{\K}{\mathcal{K}}

\begin{document}
\title[MZ inequalities and interpolation]{Marcinkiewicz-Zygmund inequalities and interpolation by spherical harmonics.}
\author {Jordi Marzo}
\address{Dept. matem\`atica aplicada i an\`alisi, Universitat de Barcelona,
                               Gran via 585, 08071 Barcelona, Spain}
\email{jmarzo@mat.ub.es}

\thanks{The author is supported by the DGICYT grant: BES-2003-2618, MTM2005-08984-C02-02
.}

\subjclass[2000]{Primary 94A20; Secondary 33C55, 65T40, 11R45}

\keywords{Spherical harmonics; Marcinkiewicz-Zygmund inequalities;
    Interpolation; Landau densities;
    Ball multiplier; Paley-Wiener spaces}

\date{\today}
\begin{abstract}
    We find necessary density conditions for
    Marcinkiewicz-Zygmund inequalities and interpolation
    for spaces of
    spherical harmonics in $\S^{d}$ with respect to the $L^{p}$ norm. Moreover, we
    prove that there are no complete interpolation families for $p\neq 2.$
\end{abstract}

\maketitle

\section{Introduction}                                      \label{SecIntro}
%--------------------------------------------------------------------------------------

    Let $\S^{d}$ be the unit sphere in $\R^{d+1}.$
    We denote by $d(u,v)=\arccos \la u, v\ra$ the geodesic distance between $u,v\in \S^{d},$
    where $\la u, v\ra$ is the scalar product in $\R^{d+1}.$
    The ball $B(u,\theta)\subset \S^{d}$
    is, therefore,
    the spherical cap
    of radius
    $0<\theta<\pi$ and center $u\in \S^{d}.$

    We consider the Banach spaces $L^{p}(\S^{d})$
    of measurable functions defined in $\S^{d}$
    such that
    $$\| f \|_{p}^{p}=\int_{\S^{d}}|f(z)|^{p}d\sigma(z)<\infty,$$
    if $1\le p< \infty,$
    and
    $$\| f \|_{\infty}=\sup_{z\in \S^{d}}|f(z)|<\infty,$$
    when $p=\infty.$ Here $d\s $ stands for the Lebesgue surface measure in $\S^{d}.$

    Now we recall some facts about spherical harmonics, see \cite{SW}.
    For any integer $\ell\ge 0,$ let $\He_{\ell}$ be the
    space of spherical harmonics of degree
    $\ell$ in $\S^{d}.$
    Then $\He_{\ell}$ is the
    restriction to $\S^{d}$
    of the homogeneous harmonic polynomials of degree $\ell$ in
    $\R^{d+1}.$

    For any integer
    $L\ge 0$ we denote the space of spherical harmonics of degree not exceeding $L$ by
    $$\Pi_{L}=\bigcup_{\ell=0}^{L}\He_{\ell}.$$
    Recall that for $p=2$ the spaces $\He_{\ell}$ are orthogonal.
    These vector spaces have dimensions
    $\dim \He_{0}=1,\dim \He_{1}=d+1$ and for all $\ell\ge 2$
    $$\dim \He_{\ell}=\frac{2\ell+d-1}{\ell+d-1}\binom{d+\ell-1}{\ell}=h_{\ell},$$
    so
    $$\dim \Pi_{L}=\frac{d+2L}{d}\binom{d+L-1}{L}=\pi_{L},$$
    and by Stirling's formula\footnote{Here and in
    what follows $\sim$ means that the ratio of the two sides is bounded from above and from
    below by two positive constants.} $\pi_{L}\sim L^{d},$ when $L\to
    \infty.$

    For any degree $L$ we take $m_{L}$ points in $\S^{d}$
    $$\mathcal{Z}(L)=\{ z_{L j}\in \S^{d}: 1\le j\le m_{L}
    \}, \;\; L\ge 0,$$
     and assume that
    $m_{L}\to \infty$ as $L\to \infty.$
    This yields a
    triangular family of points $\mathcal{Z}=\{ \mathcal{Z}(L) \}_{L\ge 0}$ in $\S^{d}.$

\begin{definition}                                                                      \label{def-MZ}
    Let $\mathcal{Z}=\{ \mathcal{Z}(L) \}_{L\ge 0}$ be a triangular
    family with $m_{L}\ge \pi_{L}$ for all $L.$
    We call
    $\mathcal{Z}$
    an $L^{p}-$\textbf{Marcinkiewicz-Zygmund} family, denoted by $L^{p}-$MZ, if there exists a constant $C_{p}>0$ such that
    for all $L\ge 0$ and $Q\in \Pi_{L},$
\begin{equation}                                                                    \label{ineq-def-MZ}
    \frac{C_{p}^{-1}}{\pi_{L}}\sum_{j=1}^{m_{L}}|Q(z_{L j})|^{p}
    \le \int_{\S^{d}}|Q(\omega)|^{p}d\sigma(\omega)
    \le \frac{C_{p}}{\pi_{L}}\sum_{j=1}^{m_{L}}|Q(z_{Lj})|^{p},
\end{equation}
    if $1\le p<\infty,$ and
    $$\sup_{\omega\in \S^{d}}|Q(\omega)|\le C\sup_{j=1,\dots , m_{L}}|Q(z_{Lj})|,$$
    when $p=\infty.$
\end{definition}

    Then the $L^{p}-$norm in $\S^{d}$ of a polynomial of degree $L$ is comparable to the discrete version
    given by the weighted $\ell^{p}-$norm of its restriction to $\mathcal{Z}(L).$
    In fact we observe that $\mathcal{Z}$ is $L^{2}-$MZ if and only
    if,
    for all $L\ge 0,$ the normalized reproducing kernels of $\Pi_{L}$
    centered at the points
    $\mathcal{Z}(L)$
    form a frame in $\Pi_{L},$ with frame bounds independent
    of $L.$

    A concept that can be seen as dual of MZ is that of interpolation.

\begin{definition}                                                                       \label{def-interp}
    Let $\mathcal{Z}=\{ \mathcal{Z}(L) \}_{L\ge 0}$ be a triangular
    family with $m_{L}\le \pi_{L}$ for all $L.$
    We say that $\mathcal{Z}$ is $L^{p}-$\textbf{interpolating}, if
    for all family $\{ c_{Lj}\}_{L\ge 0, 1\le j\le m_{L}}$ of values
    such that
    $$\sup_{L\ge 0}\frac{1}{\pi_{L}}\sum_{j=0}^{m_{L}}|c_{Lj}|^{p}<\infty,$$
    there exists a polynomial $Q\in \Pi_{L}$ such that
    $Q(z_{Lj})=c_{Lj},$ $1\le j\le m_{L}.$
\end{definition}

    Roughly speaking in order to recover the
    $L^{p}-$norm of a polynomial of degree $L$
    from the evaluation at the points in $\mathcal{Z}(L)$
    we need a sufficiently big number of points in $\mathcal{Z}(L).$
    On the other hand,
    it is possible to have a spherical harmonic of degree at most $L$
    attaining some prescribed values on $\mathcal{Z}(L)$ only when
    $\mathcal{Z}(L)$ is sparse. When we have both MZ and
    interpolation the points of the family can be thought as placed in some sort of equilibrium.

\begin{definition}                                                                  \label{CIS}
    Let $\mathcal{Z}=\{ \mathcal{Z}(L) \}_{L\ge 0}$ be a triangular family.
    We say that $\Z$ is an $L^{p}-$complete interpolating family if it is both $L^{p}-$MZ
    and $L^{p}-$interpolating.
\end{definition}

    A first measure of sparsity is the uniform separation between points of the
    same generation. This leads to the following definition.

\begin{definition}
    A triangular family $\mathcal{Z}$ is \textbf{uniformly separated} if there is a positive number $\e>0$
    such that
    $$d(z_{Lj},z_{Lk})\ge
    \frac{\epsilon}{L+1},\;\; \mbox{if}\;\; j\neq k,$$
    for all $L\ge 0.$
\end{definition}

    The precise formulation
    of the sparsity requirement is expressed in terms of
    the following Beurling type densities, \cite{OS}.

\begin{definition}
    For $\mathcal{Z}$ a triangular family in $\S^{d}$ we define the upper and lower
    density respectively as
    $$D^{-}(\Z)=\liminf_{\a\to \infty}
    \liminf_{L\to \infty}
    \frac{\min_{z\in \S^{d}}\# (\Z(L) \cap B(z,\frac{\alpha}{L+1}))}{\alpha^{d}},$$
    $$D^{+}(\Z)=\limsup_{\a\to \infty}
    \limsup_{L\to \infty}
    \frac{\max_{z\in \S^{d}}\# (\Z(L) \cap B(z,\frac{\alpha}{L+1}))}{\alpha^{d}}.$$
\end{definition}

    Now we can formulate our main result which we will prove in
    section~\ref{SecProofs}.

\begin{theorem}                                                                     \label{teo}
    Let $1\le p\le \infty.$ If $\mathcal{Z}$ is an $L^{p}-$Marcinkiewicz-Zygmund family
    there exists a
    uniformly separated $L^{p}-$MZ
    family $\tilde{\mathcal{Z}}\subset \mathcal{Z}$
    such that
    $$D^{-}(\tilde{\mathcal{Z}})\ge
    \frac{2}{d!d\sqrt{\pi}}\frac{\Gamma(\frac{d+1}{2})}{\Gamma(\frac{d}{2})}.$$
    If $\mathcal{Z}$ is an $L^{p}-$interpolating family then it is
    uniformly separated and
    $$D^{+}(\mathcal{Z}) \le
    \frac{2}{d!d\sqrt{\pi}}\frac{\Gamma(\frac{d+1}{2})}{\Gamma(\frac{d}{2})}.$$
\end{theorem}

    This result together with Theorem \ref{interp-unif}, that shows
    that a interpolating family has to be uniformly separated,
    proves that
    $L^{p}-$complete interpolation families
    must have
    $$D^{-}(\mathcal{Z})=D^{+}(\mathcal{Z})=
        \frac{2}{d!d\sqrt{\pi}}\frac{\Gamma(\frac{d+1}{2})}{\Gamma(\frac{d}{2})}.$$

    In order to stress the relationship between our problem and the
    problems of sampling and interpolation in the
    Paley-Wiener space, $PW^{p},$ of
    $L^{p}-$functions bandlimited to the unit ball, we recall some
    results. A reference for material on sampling and interpolation
    is \cite{Se}.

    As in the Paley-Wiener case, in the study of $L^{p}-$MZ and interpolation families
    much more is known
    in $d=1$ than in $d>1.$
    The main reason for such gap is that for $d=1$
    the family given by the roots of the unity is both MZ and
    interpolating.
    We recall the
    classical result due to A. Zygmund and J.
    Mar\-cin\-kiewicz: there exists a constant $C_{p}>0$ such that for any $q$ trigonometric polynomial of
    degree smaller or equal than $n$
    $$\frac{C_{p}^{-1}}{n}\sum_{j=0}^{n}|q(\o_{n,j})|^{p}\le \int_{0}^{2\pi}|q(e^{i\theta})|^{p}d\theta
    \le \frac{C_{p}}{n}\sum_{j=0}^{n}|q(\o_{n,j})|^{p},$$
    where $\o_{n,j}$ are the $(n+1)$-th roots of the unity,
    see \cite[Theorem 7.5, Chapter X]{zyg}.

    In the case $d>1$ that we deal with in this paper
    we don't have an
    even distribution of points analogous to the roots of unity, although a lot of schemes
    have been proposed. We refer to N. J. A. Sloane \cite{Slo} for further information.
    In fact, in contrast with the situation for $d=1$
    we will prove the following result about complete interpolating families.

\begin{theorem}                                                 \label{noCIS}
    For $d>1,$ there are no $L^{p}-$complete interpolating families if $p\neq 2.$
\end{theorem}

    The one dimensional case was treated by A. Zygmund and J. Mar\-cin\-kiewicz
    and can be seen as the $\S^{1}$ analogue to the Whittaker-Kotelnikov-Shannon theorem.
    Moreover,
    there is a
    complete characterization for $L^{p}-$complete interpolating families in terms of Muckenhoupt's condition,
    due to C. K. Chui, X-C. Shen and L. Zhong \cite{CSZ,CZ} analogous to that of B. S. Pavlov,
    Y. I. Lyubarskii and K. Seip \cite{Pa,LS}, in the case of the Paley-Wiener space.

    Also the classical results, for $d=1,$ about sampling and interpolation
    for Bernstein's space given by A. Beurling \cite{Beu}
    using densities and weak limits
    have their counterparts for $L^{p}-$MZ and interpolation families in the recent results given in \cite{OS}.
    Indeed, it is shown in \cite{OS} that if a triangular family is
    $L^{p}-$MZ then its lower density has to be greater
    or equal to $1/2\pi,$ and that the converse holds for families with densities
    greater to $1/2\pi.$
    The corresponding result for interpolation can be proved without
    a lot of effort.

    In the Paley-Wiener case and for greater dimensions
    there are classical necessary
    conditions for sampling and interpolation in terms of densities due to
    H. Landau \cite{LAN}. It can be easily seen that
    these densities can not characterize sampling and interpolation
    sequences.
    In previous work \cite{Jo} we have shown how to obtain
    sampling and interpolation sequences with densities arbitrarily close to the critical
    one ( Nyquist density) for functions bandlimited in the
    Euclidean space. In particular this applies to functions
    in $PW^{p}.$

    Concerning the question of sufficient conditions in $\S^{d},$
    in 2000 H. N. Mhaskar, F. J. Narcowich and J. D. Ward \cite{MNW}
    using the doubling weights construction due to G. Mastroianni and V. Totik \cite{MT}
    obtained  a sufficient condition for being
    $L^{p}-$MZ in terms of a mesh norm condition that is far from being optimal.

    Our main
    result,  Theorem~\ref{teo}, can be seen as the analogue of the Paley-Wiener space result due to H. Landau \cite{LAN}.
    Instead of using the
    approach provided by
    J. Ramanathan T. Steger \cite{RS}, that was adapted in \cite{OS}
    to the $\S^{1}$ case, we are going to adapt
    the classical operator theoretic proof given by H. Landau.
    We deal with the case $d>1$ but the
    result for $d=1$ follows also with minor
    changes.

    We prove also that for $p\neq 2$ there are no triangular
    families that are both $L^{p}-$MZ and interpolating.
    Indeed, if such a family exists one can construct a
    bounded multiplier that turns out to be the multiplier for the
    ball. Finally the well known result of C. Fefferman \cite{F} brings us the
    contradiction.

    Up to here we have seen that the knowledge is similar in both spaces.
    Therefore the Paley-Wiener case provides us the inspiration but
    technically the situation is completely different.
    In further work we will focus on this relation.

    The main technical difficulties in the case $d>1$
    is that
    we can't use the techniques for holomorphic polynomials used in
    \cite{OS}, like Hadamard's three sphere theorem, Bernstein type
    inequality or sub-mean value inequality.

    The outline of this paper is as follows. In the next section we summarize some well known
    facts about spherical
    harmonics and Jacobi polynomials.

    In section \ref{SecConcentr} we calculate
    the trace of the concentration operator over a spherical cap and his power, which are the main tools
    in proving the density conditions.
    Controlling  these quantities
    we can estimate the number of "big" eigenvalues of the
    concentration operator, and this quantity can be thought of as
    the local dimension of the space of spherical harmonics. Now, to
    get a MZ or interpolating family we will need locally to have
    either more or less points
    than this local dimension.

    In section \ref{SecFam} we prove several general results
    concerning MZ and interpolating families.
    Our main tool, Lemma \ref{comparacio}, says that the $L^{p}-$norm
    of a spherical harmonic is equivalent, with constants that do not depend on the degree, to the
    $L^{p}-$norm computed in any other sphere with radius close to
    1.
    A perturbative argument
    allow us to treat only the case $p=2$
    with uniformly separated family in Theorem \ref{teo}.
    We characterize also
    the Carleson families of measures in $\S^{d}.$

    In section \ref{secnoCIS} we prove the result about nonexistence
    of complete interpolating families,
    Theorem \ref{noCIS},
    using the approach outlined above.

    Finally, in section \ref{SecProofs} we prove two technical lemmas
    that we use to prove
    the main result.

\section{Spherical Harmonics}                                      \label{SecSpheHarm}
%--------------------------------------------------------------------------------------

    In this section we recall some facts about spherical harmonics
    and Jacobi polynomials, see \cite{SW,S}.

    Let $Z^{\ell}_{\eta}\in \He_{\ell}$ be
    such that for
    $Q\in \He_{\ell}$
    $$Q(\eta)=\int_{\S^{d}}Q(\xi)\overline{Z^{\ell}_{\eta}(\xi)}d\s (\xi),\;\;\; \xi \in \S^{d}.$$
    We call it the {\it zonal harmonic} of degree $\ell$ with pole $\eta\in \S^{d}.$
    Let $\mathcal{P}(\S^{d})$ be the linear span of $\bigcup_{\ell=0}^{\infty}\He_{\ell}.$

\begin{definition}
    We call \textbf{zonal multiplier} any linear map from $\mathcal{P}(\S^{d})$ into $\mathcal{C}(\S^{d})$
    which commutes with rotations.
\end{definition}

    The following explains why the term multiplier is used in
    this last definition.

\begin{theorem}{\cite[chap. 3]{CW}}                                                 \label{co_we}
    Let $T$ be a zonal multiplier in $\S^{d}.$ For any $\ell\ge 0,$
    $Y_{\ell}\in \He_{\ell}$ are eigenvectors of $T$ corresponding to the same eigenvalue.
\end{theorem}

    Then for $T$ as above there exists a sequence $\{ m_{\ell} \}_{\ell=0}^{\infty}\subset \C$
    such that for $\sum_{\ell=0}^{N}Y_{\ell}  \in \mathcal{P}(\S^{d})$
    $$T\left(\sum_{\ell=0}^{N}Y_{\ell}\right)=\sum_{\ell=0}^{N}m_{\ell}Y_{\ell}.$$

\begin{definition}
    We say that $T$ is a \textbf{bounded zonal multiplier}
    if for some $1\le p<\infty$ we have $A_{p}>0$ such that for any
    $Y\in \mathcal{P}(\S^{d})$
    $$\| TY \|_{p}\le A_{p} \| Y \|_{p}.$$
\end{definition}

\begin{definition}
    We call a function in $\S^{d}$ \textbf{zonal} if it is invariant by the action of $SO(d),$
    i.e. if
    $$f\circ \rho (\o)=f(\o),\;\;\; \o\in \S^{d},$$
    for $\rho\in SO(d+1)$ such that $\rho N=N.$
\end{definition}

    Observe that this is equivalent to saying that $f$ is constant on
    $$L_{\theta}=\{\o\in \S^{d}: d(\o,N)=\theta \},\;\;\; 0\le \theta\le \pi,$$
    so the value of a zonal function in one point depends only on its geodesic distance
    to the north pole.

    For functions $f,g\in L^{1}(\S^{d})$ with $g$ zonal
    we define the convolution product
    $$(g\ast f)(\o)=\int_{\S^{d}}g^{\flat}(\la \o , x\ra )f(x)d\s (x),$$
    where $g^{\flat}$ is the function in $[-1,1]$ defined by
    $$g^{\flat}(\la \o, N\ra )=g(\o).$$

    In the Hilbert space $L^{2}(\S^{d})$
    we can take an orthonormal basis of $\He_{\ell},$ that we denote by
    $Y_{\ell}^{1},\dots ,Y_{\ell}^{h_{\ell}},$
    which can be chosen in such a way that $Y_{\ell}^{1}$ is the only vector non-vanishing at the north
    pole.
    Taking all these basis for $\ell=0,\dots L$ together we get an orthonormal basis for $\Pi_{L}.$
    Given $f\in L^{2}(\S^{d})$ we define its Fourier coefficients as
    the triangular family
    $$\hat{f}(\ell,j)=\int_{\S^{d}}f(z)\overline{Y_{\ell}^{j}(z)}d\s (z),$$
    for $\ell\ge 0 $ and $1\le j\le h_{\ell}.$

    It is well known that the reproducing kernel for $\Pi_{L}$ is
    $$K_{L}(u,v)=\sum_{\ell=0}^{L}\sum_{j=1}^{h_{\ell}}Y_{\ell}^{j}(u)\overline{Y_{\ell}^{j}(v)},\;\;\;u,v\in \S^{d},$$
    and that this expression does not depend on the basis.

    Now we will compute the kernel $K_{L}.$
    The zonal harmonic of degree $\ell\ge 0$ is the reproducing
    kernel in $\mathcal{H}_{\ell},$ so
    $$Z^{\ell}_{u}(v)=\sum_{j=1}^{h_{\ell}}Y_{\ell}^{j}(u)\overline{Y_{\ell}^{j}(v)}=\frac{h_{\ell}}{\s(\S^{d})}
    P_{\ell}(d+1;\la u,v\ra ),$$ where $P_{\ell}(d+1;x)$ is the $\ell-$th Legendre polynomial
    in $d+1$ dimensions, \cite{Mu}.
    Using the Christoffel-Darboux formula we get
\begin{align*}
    \sum_{\ell=0}^{L} & \frac{h_{\ell}}{\s(\S^{d})}P_{\ell}(d+1;\la u, v\ra )
    =
    \binom{d+L-1}{L}
    \frac{P_{L}(d+1;\la u, v\ra)-P_{L+1}(d+1;\la u, v\ra)}{\s(\S^{d})(1-\la u, v\ra)}.
\end{align*}
    Finally,
    $$P_{L} (d+1;x)-P_{L+1}(d+1;x)=(1-x)\binom{L+(d-2)/2}{L}^{-1}P_{L}^{(d/2,(d-2)/2)}(x),$$
    where
    $P_{L}^{(\alpha,\beta)}$ stands for the Jacobi polynomial of degree $L$ and
    index $(\alpha,\beta).$

    From now on we denote $\l=(d-2)/2.$
    So the reproducing kernel is given by
    $$K_{L}(u,v)=\frac{C_{d,L}}{\s(\S^{d})}P_{L}^{(1+\l,\l)}
    (\la u, v\ra),$$
    where
    $C_{d,L}=\binom{d+L-1}{L}/\binom{L+\frac{d-2}{2}}{L},$
    and using Stirling's formula one can see that $C_{d,L}\sim L^{d/2},$ if $L\to \infty.$

    To
    estimate the $L^{p}$ norm of this kernel, all we need is to estimate the
    $L^{p}-$norm of the Jacobi polynomial.
    For the case $p=\infty$ it is well known that
    $$\sup_{t\in [-1,1]} | P_{L}^{(1+\l,\l)}(t)|= \binom{L+\l+1}{L}
    \sim L^{d/2}.$$
    For $1\le p<\infty$ we can use the estimate in \cite[p. 391]{S} and the
    fact that
    $P_{L}^{(1+\l,\l)}(t)=(-1)^{L}P_{L}^{(1+\l,\l)}(-t)$
    to obtain, for any $v\in \S^{d}$

\begin{equation}                                                                        \label{norma}
\int_{\S^{d}}|P_{L}^{(1+\l,\l)}(\la u, v\ra)|^{p}d\sigma(u)\sim
\left\{ \begin{array}{ll}
    L^{d(\frac{p}{2}-1)},
    & p>\frac{2d}{d-1},
    \\
    \\
    L^{-p/2}\log L,
    & p=\frac{2d}{d+1},
    \\
    \\
    L^{-p/2},
    &
    p<\frac{2d}{d+1}.
     \end{array}\right.
\end{equation}

    Finally we recall an estimate that will be used later on
    \cite[p. 198]{S}:
\begin{equation}                                                                        \label{estimate}
    P_{L}^{(1+\l,\l)}(\cos
    \theta)=\frac{k(\theta)}{\sqrt{L}}\left\{ \cos \left((L+\l+1)\theta+\gamma \right)+\frac{O(1)}{L\sin
    \theta}\right\},
\end{equation}
    if $c/L\le \theta \le \pi-(c/L),$ where
    $$k(\theta)=\pi^{-1/2} \left( \sin \frac{\theta}{2}\right)^{-\l-3/2}\left( \cos
    \frac{\theta}{2}\right)^{-\l-1/2},\;\; \gamma=-\left(\l+\frac{3}{2}\right)\frac{\pi}{2}.$$

\section{Concentration Operator}                                      \label{SecConcentr}
%--------------------------------------------------------------------------------------

    In this section we estimate the trace of
    the concentration operator and its square in order to obtain an estimate for the eigenvalues of
    this operator, Proposition
    \ref{prop}.
    In the next section we will show how the cardinality of the set of "big" eigenvalues
    can be related with the density of the triangular family when it is MZ
    or interpolating.

    Let $\K_{A}$ be the concentration operator over
    $A\subset\S^{d}$ defined for $Q\in \Pi_{L}$ and given by
\begin{equation}                                                                \label{concentr}
    \K_{A} Q(u)=\int_{A}K_{L}(u,v)
    Q(v)dv.
\end{equation}
    This operator
    results from the composition of the restriction operator
    $$\begin{array}{ccc} \Pi_{L} & \longrightarrow & L^{2}(\S^{d}) \\
Q & \longmapsto & \chi_{A}Q,
\end{array}
$$
    with the orthogonal projection
    $$\begin{array}{ccl}L^{2}(\S^{d}) & \longrightarrow & \Pi_{L} \\
f & \longmapsto & \sum_{\ell=0}^{L}\sum_{j=1}^{h_{\ell}}\la f,
Y_{\ell}^{j}\ra Y_{\ell}^{j}.
\end{array}
$$
    The operator $\K_{A}$ is self-adjoint and by the spectral theorem
    its eigenvalues are all real and $\Pi_{L}$
    has an orthonormal basis of eigenvectors of $\K_{A}$.
    We can compute the trace of this operator using $Z_{u}^{\ell}(u)=h_{\ell}/\s(\S^{d})$
    and the expression of $K_{L}$ as sum of zonal harmonics
    $${\rm tr} (\K_{A})=\int_{A}K_{L}(u,u)
    d\s(u)=\pi_{L}\frac{\s(A)}{\s(\S^{d})}.$$

    Now we take $A$ a spherical cap with radius $\alpha/(L+1)$
    and we want to obtain an estimate for ${\rm tr} (\K_{A}^{2}).$

\begin{proposition}                                                                 \label{prop}
    Let $A\subset \S^{d}$ be a spherical cap with radius $\a/(L+1)$
    and let $\K_{A}$ be the concentration operator defined in
    (\ref{concentr}). Then
    $${\rm tr} (\K_{A})-{\rm tr} (\K^{2}_{A})=O(\a^{d-1}\log \a),$$
    when $L\to \infty,$ with constants depending only on $d.$
\end{proposition}

\remark
    The invariance of the zonal harmonic, $Z_{\rho u}^{\ell}(\rho v)=Z_{u}^{\ell}(v),$ for $\rho\in SO(d+1),$
    gives
    ${\rm tr} (\K^{2}_{A})={\rm tr} (\K_{\rho A}^{2}).$

\proof
    Using the reproducing property we have
\begin{align*}
    {\rm tr} & (\K^{2}_{A})=\int_{A}\int_{A}|K_{L}(u,v)|^{2}d\s(u)d\s(v)
    \\
    &
    =
    \int_{A}\int_{S^{d}}|K_{L}(u,v)|^{2}d\s(u)d\s(v)
    -\int_{A}\int_{S^{d}\setminus A}|K_{L}(u,v)|^{2}d\s(u)d\s(v)
    \\
    &
    =
    \int_{A} K_{L}(u,u) d\s(u)
    -\int_{A}\int_{S^{d}\setminus A}|K_{L}(u,v)|^{2}d\s(u)d\s(v)
    \\
    &
    =
    {\rm tr}(\K_{A}) -\frac{C_{d,L}^{2}}{\s(\S^{d})^{2}}\int_{A}\int_{\S^{d}\setminus A}
    | P_{L}^{(1+\l,\l)}(\la u, v\ra )|^{2}d\s(u) d\s(v).
\end{align*}
     In $\S^{d}$ we take the spherical coordinates
    $$\left\{ \begin{array}{ll}
    x_{1}=\sin \theta_{d}\dots \sin \theta_{2}\sin \theta_{1},
    \\
    x_{2}=\sin \theta_{d}\dots \sin \theta_{2}\cos \theta_{1},
    \\
    \dots
    \\
    x_{d}=\sin \theta_{d}\cos \theta_{d-1},
    \\
    x_{d+1}=\cos \theta_{d},
\end{array}\right. $$
    where $0\le \theta_{k}<\pi$ if $k\neq 1$ and $0\le \theta_{1}<2\pi.$
     Using the rotation invariance
     we get
\begin{align*}
    \int_{\S^{d}\setminus A}
     | P_{L}^{(1+\l,\l)}(\la u, v\ra)|^{2}d\s(u)\le & \int_{\S^{d}\setminus B(N,d(v,\partial A))}
     | P_{L}^{(1+\l,\l)}(\la u, N\ra)|^{2}d\s(u)
     \\
     =
     \s(\S^{d-1}) &
    \int_{d(v,\partial A)}^{\pi}
     |P^{(1+\l,\l)}_{L}(\cos \theta)|^{2}\sin^{d-1}\theta \,d\theta.
\end{align*}
    Let $\theta_{\a}=\a/(L+1)$ be the radius of the spherical cap $A$
    and
    let $v\in A$ be fixed. Since we want an asymptotic result we will take an $\a\gg
    1$ and an even bigger $L,$ in such a way that $\theta_{\a}<<1.$
    Integrating over $A$ we get
\begin{align*}
    \int_{A} & \int_{\S^{d}\setminus A}
     | P_{L}^{(1+\l,\l)}(\la u, v\ra )|^{2}d\s(u)d\s(v)
     \\
     \le
     &
    \s(\S^{d-1})^{2}
     \int_{0}^{\theta_{\a}}\sin^{d-1}\eta
    \int_{\theta_{\a}-\eta}^{\pi}
     |P^{(1+\l,\l)}_{L}(\cos \theta)|^{2}\sin^{d-1}\theta \,d\theta
     d\eta.
\end{align*}
    Split the innermost integral depending on whether $\theta>L^{-1}$ or $\theta <L^{-1}.$
    In the first case (obs. $\theta_{\a}>L^{-1}$)
\begin{align*}
    L^{d}\int_{0}^{\theta_{\a}} & \sin^{d-1} \eta
    \int_{\theta_{\a}-\eta,\theta>L^{-1}}^{\pi}
    |P_{L}^{(1+\l,\l)}(\cos \theta)|^{2}\sin^{d-1} \theta \,d\theta d\eta
    \\
    &
    \lesssim \int_{0}^{\a} \eta^{d-1}
    \int_{\pi-m(\a,\eta,L)}^{\pi}
    |P_{L}^{(1+\l,\l)}(\cos \theta)|^{2}\sin^{d-1} \theta \,d\theta d\eta
    \\
    &
    +
    \int_{0}^{\alpha} \eta^{d-1}
    \int_{m(\a,\eta,L)}^{\pi-m(\a,\eta,L)}
    |P_{L}^{(1+\l,\l)}(\cos \theta)|^{2}\sin^{d-1} \theta \,d\theta
    d\eta
    =A1+A2,
\end{align*}
    where $m(\a,\eta,L)=\max((\a-\eta)/L,1/L).$

    For part $A1$ we use that
    $|P^{(1+\l,\l)}_{L}(x)|=O(L^{\l}),$
    for $-1\le x\le 0,$ \cite[p. 168]{S}.
    Then, for a fixed $\a$
    $$A1\lesssim L^{2\l}\int_{0}^{\a}\eta^{d-1}m(\a,\eta,L)^{d}
    d\eta=
    L^{-2}\int_{0}^{\a}\eta^{d-1}\max(\a-\eta,1)^{d}  d\eta,$$
    which goes to zero as $L\to \infty.$

    Using the Szeg\"o estimate (\ref{estimate}) we get
\begin{align*}
    A2 & \lesssim \int_{0}^{\a}\eta^{d-1} \int_{
    m(\alpha,\eta,L)}^{\pi-
    m(\alpha,\eta,L)}\frac{k^{2}(\theta)}{L}\sin^{d-1} \theta d\theta
    d\eta
    =
    \int_{0}^{\a}\eta^{d-1} \int_{
    m(\alpha,\eta,L)}^{\pi-
    m(\alpha,\eta,L)}\frac{2^{d-1}}{L\sin^{2} \frac{\theta}{2}} d\theta
    d\eta
    \\
    &
    \sim
    \frac{1}{L}\int_{0}^{\a}\eta^{d-1} \cot
    \frac{m(\alpha,\eta,L)}{2}d\eta
    \lesssim \frac{1}{L}\cot
    \frac{1}{L}\int_{\a-1}^{\a}\eta^{d-1} d\eta
    +
    \int_{1}^{\a}(\a-\eta)^{d-1}\frac{1}{L}\cot \frac{\eta}{L}d\eta
    \\
    &
    \lesssim
    \frac{\a^{d-1}}{L}\cot
    \frac{1}{L}
    +\int_{1}^{\a}\frac{(\a-\eta)^{d-1}}{\eta}d\eta= O(\a^{d-1} \log \a).
\end{align*}
\\
    For the second part ($\theta<L^{-1}$) we obtain
\begin{align*}
    \int_{0}^{\theta_{\a}} & \sin^{d-1} \eta \int_{\theta_{\a}-\eta,
    \theta<L^{-1}}^{\pi}
    |P_{L}^{(1+\l,\l)}(\cos \theta)|^{2}\sin^{d-1} \theta \,d\theta d\eta
    \\
    &
    =
    \int_{\theta_{\a}-L^{-1}}^{\theta_{\a}} \sin^{d-1} \eta \int_{\theta_{\a}-\eta}^{L^{-1}}
    |P_{L}^{(1+\l,\l)}(\cos \theta)|^{2}\sin^{d-1} \theta \,d\theta
    d\eta.
\end{align*}
    Observe that $\eta<\theta_{\a}-L^{-1}$ would imply $\theta>L^{-1}.$ Then
\begin{align*}
    L^{d} & \int_{\theta_{\a}-L^{-1}}^{\theta_{\a}} \sin^{d-1} \eta \int_{\theta_{\a}-\eta}^{L^{-1}}
    |P_{L}^{(1+\l,\l)}(\cos \theta)|^{2}\sin^{d-1} \theta \,d\theta d\eta
    \\
    &
    \le
    L^{2d}
    \int_{\theta_{\a}-L^{-1}}^{\theta_{\a}}
    \sin^{d-1} \eta \int_{\theta_{\a}-\eta}^{L^{-1}}\sin^{d-1} \theta \,d\theta d\eta
    \sim \int_{\alpha-1}^{\a}(1-(\a-t)^{d})t^{d-1}dt=O(\a^{d-1}).
\end{align*}
    Taking all the estimates together we get the result.
\qed

\section{General results about MZ and interpolating families}                    \label{SecFam}
%--------------------------------------------------------------------------------------

    In this section we prove some results about
    MZ and interpolation triangular families. Also we characterize the
    families of
    Carleson measures
    for the spherical harmonics $\Pi_{L}$
    on $\S^{d}.$

    The first thing we need to show is that in calculating densities
    we can restrict
    ourselves to uniformly
    separated families.
    Following \cite{OS}
    we will compare the norm of a polynomial in $\S^{d}$ with the norm in a
    shell sufficiently small containing $\S^{d}.$
    This comparison
    is harder than in dimension one \cite[Lemma 2]{OS} because
    Hadamard's three circle principle is no longer available.

    For $r>0$ we denote $S^{d}_{r}=r\S^{d}$ and for a measurable function $f$ defined in $S^{d}_{r}$
    we have
    $$\frac{1}{r^{d}}\int_{S^{d}_{r}}f(\omega)d\s(\omega)=\int_{S^{d}}f(r\omega)d\s(\omega).$$
    First we prove a result which we will use later one.

\begin{proposition}                                                                         \label{multi}
    There exists a bounded zonal multiplier $T:L^{p}(\S^{d})\longrightarrow L^{p}(\S^{d})$
    for $1\le p\le\infty,$
    such that
    $\| T \|_{p}\le C<\infty,$ with $C$ independent of $p$ and $L,$
    and such that $\rank \, T \subset \Pi_{3L},$ $T |\Pi_{L}= Id.$
\end{proposition}

\proof
    Let $g\in L^{1}(\S^{d})$ be a zonal function. For any $1\le p\le \infty$ we have
    $$\| g\ast f \|_{p}\le \| g \|_{1} \| f \|_{p},$$
    so the operator $T_{g}:\mathcal{P}(\S^{d})\longrightarrow \mathcal{C}(\S^{d})$ defined as
    $T_{g}(f)=g\ast f$ is bounded in $L^{p}(\S^{d}),$ commutes with rotations and has norm
    $\| g \|_{1}.$

    Using H\"older's inequality it is easy to see that the function
    $$g=\frac{\binom{2L+\l+1}{2L}}{\binom{L+\l+1}{L}} P_{L}^{(1+\l,\l)}(\la N, \cdot \ra)
    P_{2L}^{(1+\l,\l)}(\la N, \cdot \ra),$$ has $L^{1}-$norm independent of $L$.
    Also, for $f\in \mathcal{P}(\S^{d})$
\begin{align*}
    g\ast f(\o) & =\int_{\S^{d}}g^{\flat}(\la \o, x\ra)f(x)d\s (x)
    \\
    &
    =
    \frac{\binom{2L+\l+1}{2L}}{\binom{L+\l+1}{L}}
    \int_{\S^{d}}P_{L}^{(1+\l,\l)}(\la \o, x\ra)P_{2L}^{(1+\l,\l)}
    (\la \o, x\ra)f(x)d\s (x),
\end{align*}
    is a polynomial of degree $\le 3L$ in $\o,$ hence $\rank \, T_{g}\subset \Pi_{3L}.$
    Finally taking the polynomial
    $f\in \Pi_{L}$ and applying the reproducing property we obtain
    $g\ast f(\o)=f(\o),$
    so $T_{g}|\Pi_{L}= Id.$
\qed
%==========================================================================================================

    The next lemma shows that the $L^{p}-$norm of a spherical harmonic in the unit sphere is equivalent
    to the $L^{p}-$norm in any other sphere with radius close to 1.

\begin{lemma}                                                                            \label{comparacio}
    Let $p\in [1,\infty]$ and $Q\in \Pi_{L}.$ For any $|r-1|\le \rho/L$
    there exists a constant $C$ depending only on $p$ and $\rho$ such
    that
\begin{equation}                                                                    \label{desigualtat}
    C \|Q\|_{L^{p}(\S^{d})}\le \| Q
    \|_{L^{p}(S^{d}_{r})} \le C^{-1} \| Q \|_{L^{p}(\S^{d})}.
\end{equation}
\end{lemma}

\proof
    First we consider the right hand side inequality.
    For $Q\in \Pi_{L},$ $|Q|^{p}$ is subharmonic, thus
    for
    $0<r<1$ and $1 \le
    p<\infty$
    $\| Q \|_{L^{p}(S^{d}_{r})}\le \| Q \|_{L^{p}(\S^{d})},$ \cite[Theorem 2.12]{HK}.
    For $p=\infty$ the same inequality follows using the
    maximum principle.

    Using the orthogonal decomposition in spherical harmonics of a harmonic functions in $\S^{d}$
    it can be proved that Hadamard's three circle principle for harmonic functions holds
    in $L^{2}-$norm
    \cite[lemma 2.1]{KM}. Then, for $Q\in \Pi_{L},$ $1<r<1+\rho/L$ and $R\gg 1,$ we have
    $$\log \| Q \|_{L^{2}(S^{d}_{r})}\le \left( 1-\frac{\log r}{\log
    R}\right)\log \| Q \|_{L^{2}(\S^{d})}+\frac{\log r}{\log
    R}\log \| Q \|_{L^{2}(S^{d}_{R})},$$
    and using that $\| Q \|_{L^{2}(S^{d}_{R})}^{2}=O(R^{L})$
    we obtain
    $\| Q \|_{L^{2}(S^{d}_{r})}\le e^{\rho}\| Q
    \|_{L^{2}(S^{d})}.$

    Let $Q_{L}\in \Pi_{L}$
    be such that $\| Q_{L} \|_{\infty}=1=Q_{L}(N)$ and let $1-\rho/L<r<1.$
    Restricting $Q_{L}$ to a great circle of $\S^{d}$ through $N$ we get a trigonometric polynomial of degree
    at most $L.$ So using Bernstein's inequality we get $Q_{L}(z) \ge 1-\e$ for all $z\in B(N,\e/L).$

    We want to estimate the integral
    $$Q_{L}(rN) =\frac{1}{\s(\S^{d})}\int_{\S^{d}}\frac{1-r^{2}}{|rN-u|^{d+1}}Q_{L}(u)d\s(u).$$
    For any $0<\theta <1-r$ we have
    $$\frac{1-r^{2}}{1+r^{2}-2r\cos \theta}\ge \frac{1}{1-r},$$
    then the integral over $B(N,\e/L)$ is bounded below by a constant independent of $r$
\begin{align*}
    \frac{1}{\s(\S^{d})} & \int_{B(N,\e/L)}\frac{1-r^{2}}{|rN-u|^{d+1}}Q_{L}(u)d\s(u)
    \\
    &
    =
    (1-\e)\frac{\s(\S^{d-1})}{\s(\S^{d})}\int_{0}^{\e/L}\left[ \frac{(1-r^{2})^{2}}{1+r^{2}-2r\cos \theta}
    \right]^{(d+1)/2}\frac{\sin^{d-1} \theta}{(1-r^{2})^{d}}d\theta
    \\
    &
    \gtrsim
    \frac{1-\e}{(1-r)^{d}} \int_{0}^{\e/L} \chi_{(0,1-r)}(\theta)\sin^{d-1}\theta d\theta \gtrsim
    (1-\e)\left(\frac{\e}{\rho}\right)^{d}.
\end{align*}
    Since
    $$\frac{(1-r^{2})^{2}}{1+r^{2}-2r\cos \theta}=\frac{(1-r^{2})^{2}}{2r(1-\cos \theta)+(1-r)^{2}}
    \le \frac{2(1-r)^{2}}{1-\cos \theta},$$
    then
    $$\int_{B(N,\e/L)^{c}}\frac{1-r^{2}}{|rN-u|^{d+1}}Q_{L}(u)d\s(u)\le C(1-r)\frac{L}{\e}\le C\frac{\rho}{\e}.$$
    We have seen that there exists a constant $\delta_{d}>0,$ depending only on $d,$
    such that for $0<\rho<\delta_{d},$ $0<1-\rho/L<r$ and $Q\in \Pi_{L}$
    $\| Q \|_{L^{\infty}(\S^{d}_{r})}\ge C_{\rho}\| Q \|_{L^{\infty}(\S^{d})}.$
    Now, iterating the process, and therefore changing the constant, we can obtain the same result for
    arbitrary $\rho>0$ getting
    for any $0<1-\rho/L<r$ and $Q\in \Pi_{L}$
    $$\| Q \|_{L^{\infty}(\S^{d}_{r})}\ge C_{\rho}\| Q \|_{L^{\infty}(\S^{d})}.$$

    So the dilation operator $T_{r}$ in $\Pi_{L}$
    given by $Q\mapsto Q(r\cdot )$ is such that, if we denote by $|T_{r} |_{p}$
    the norm of $T_{r}$ defined in $(\Pi_{L},\| \cdot \|_{p}),$ we get
    $| T_{r} |_{2}\le e^{\rho}$
    and $| T_{r} |_{\infty}\le C_{\rho}.$
    Being $\Pi_{L}$ finite dimensional spaces we always have
    $| T_{r} |_{p}<\infty.$ By \cite[Theorem
    VI.10.10,p.524]{DS} we know that
    $\log | T_{r} |_{p}$ is a convex function of $1/p,$ then
    for all $2\le p\le \infty$ we have
    $| T_{r} |_{p}\le \max\{C_{\rho},e^{\rho}\}.$

    For $1< p<2$ we consider the
    multiplier $M=M_{L}$ given by Proposition \ref{multi}.
    Then for $Q\in \Pi_{L}$ and $1<p<2$
\begin{eqnarray*}
    \| T_{r}(Q) \|_{L^{p}(\S^{d})} & = & \sup_{\| R
    \|_{q}\le 1}\left| \la T_{r}(Q),R \ra  \right|
    =
    \sup_{\| R
    \|_{q}\le 1}\left| \sum_{k=0}^{L}r^{k}\la M_{L}Q_{k},R \ra \right|
    \\
    & = &
    \sup_{\| R
    \|_{q}\le 1}\left| \sum_{k=0}^{L}r^{k}\la Q_{k},M_{L}R \ra \right|
    \lesssim
    \sup_{\| R
    \|_{q}=1}\| Q \|_{p} \| T_{r}(M_{L} R) \|_{q}
     \\
    & \lesssim &
    | T_{r}|_{q} \| Q \|_{p}\le C_{\rho} \| Q \|_{p}.
\end{eqnarray*}
    We observe that we can't use the projection onto $\Pi_{L}$
    instead of $M_{L}$ in the calculation above
    because for $p\neq 2$ it is not bounded by a constant independent of $L,$ see section \ref{secnoCIS}.

    So far we have seen that for $1<p\le\infty,$ $1<r<1+\rho/L$ and $Q\in \Pi_{L}$
    $$\|Q \|_{L^{p}(\S^{d})}\le C_{\rho}\| Q \|_{L^{p}(\S^{d})}.$$
    For $p=1$ we can just take the limit.

    For the left hand side inequality in (\ref{desigualtat}) with $r>1$
    we define, given $Q\in \Pi_{L},$ the polynomial $\tilde{Q}(\omega)=Q(r \omega)$
    and apply the former result.
\qed
%==========================================================================================================

    Integrating with respect to the radius we get the following
    analog of \cite[Corollary1]{OS}.

\begin{corollary}                                                                   \label{corl}
    Let
    $$C_{\rho,L}=\{ \omega \in \R^{d+1}:||\omega|-1|<\rho/L \}.$$
    For $Q\in \Pi_{L}$ and $1\le  p\le  \infty$ we have
    $$\| Q \|_{L^{p}(\S^{d})}^{p}\asymp L\| Q \|_{L^{p}(C_{\rho,L},dm)}^{p},$$
    where the constants depend on $\rho$ and $p,$ but not on the polynomial.
\end{corollary}

    Now we want to prove that a triangular family $\mathcal{Z}$ is
    uniformly separated if and only if the left hand inequality in (\ref{ineq-def-MZ}) holds.
    This is the generalization to $d\ge 1$ of \cite[Theorem 3]{OS}
    and will be used to show
    that a $MZ$ family contains a separated family which is also $MZ.$
    The problem in proving this result comes from the fact that
    there is no analogue of the
    Bernstein inequality
    for spherical harmonics if $p\neq \infty.$
    Instead of proving our result directly,
    we will derive it
    from the next characterization for Carleson
    measures on $\S^{d}$ that can be of interest on their own.

\begin{definition}                                          \label{def-carleson}
    Let M=$\{ \mu_{L} \}_{L\ge 0}$ a family of measures on $\S^{d}$ and $1\le p<\infty.$ We say that $M$ is an $L^{p}$-\textbf{Carleson}
    family
    for $\Pi_{L}$
    if there exists a positive constant $C$ such that for any $Q\in \Pi_{L}$
    $$\int_{\S^{d}}|Q(z)|^{p}d\mu_{L}(z)\le C \int_{\S^{d}}|Q(z)|^{p}d\s(z).$$
\end{definition}

\begin{theorem}                                                     \label{carleson}
    Let $1\le p<\infty.$ The family of measures M=$\{ \mu_{L} \}_{L\ge 0}$ on $\S^{d}$
    is $L^{p}-$Carleson for $\Pi_{L}$ if and only if there exists a $C>0$ such that
\begin{equation}                                                            \label{cond-carleson}
    \sup_{z\in \S^{d}}\mu_{L}(B(z,L^{-1}))\le \frac{C}{\pi_{L}}.
\end{equation}
\end{theorem}

\remark We want to point out that condition (\ref{cond-carleson}) is
independent of $p$ and that we could take
    balls of any other radius $\alpha /L $ for $\a>0.$

\proof
    Let $0<m_{d}$ be the first extremum of the Bessel function $J_{d/2}$
    and let $\eta_{L}$ be such that $\eta_{L}L\to m_{d}$ when $L\to \infty.$
    Now, using Mehler-Heine formula \cite[Theorem 8.1.1.]{S} we see that
    there exist $\delta_{d}>0$ and $L_{0}$ such that for $L\ge L_{0}$
    and
    $0\le \eta\le \eta_{L}$
    $$1\ge L^{-d/2}P_{L}^{(1+\l,\l)}(\cos \eta)\ge L^{-d/2}P_{L}^{(1+\l,\l)}(\cos \eta_{L})\ge \delta_{d}>0.$$
    We argue by contradiction.
    Suppose that for all $n\in \N$ there exist $L_{n}$ and a geodesic ball $B_{n}$
    with radius $m_{d}/L_{n}$ such that $\pi_{L_{n}}\mu_{L_{n}}( B_{n} )>n.$
    Let $b_{n}\in \S^{d}$ be the center of $B_{n}$ and define for $\w\in \S^{d}$
    $$K_{n}(\w)=P_{L_{n}}^{(1+\l,\l)}(\la b_{n}, \w\ra )\in \Pi_{L_{n}}.$$
    For any Carleson family of measures $M$ we get
\begin{align*}
    \| L_{n}^{-d/2}K_{n} \|_{p}^{p} & \gtrsim \int_{\S^{d}}|L_{n}^{-d/2}K_{n}(z)|^{p}d\mu_{L_{n}} (z)
    \ge
    \int_{B_{n}}|L_{n}^{-d/2}K_{n}(z)|^{p}d\mu_{L_{n}} (z)\ge \delta_{d}^{p}\mu_{L_{n}}(B_{n}).
\end{align*}
    Then $L_{n}^{-d(p/2-1)}\| P_{L_{n}}^{(1+\l,\l)}(\la b_{n},\cdot \ra )\|_{p}^{p}\ge C n$
    with $C$ depending on $p$ and $d,$
    so if we take
    $p\ge 2d/(d+1)$ this
    contradicts (\ref{norma}).

    For other $p\ge 1$ we consider $\ell$ such that $q=\ell p>2d/(d+1).$ Then for
    $$K_{n}(\omega)=P_{[L_{n}/\ell]}^{(1+\l,\l)}(\la b_{n}, \w\ra )^{\ell}\in \Pi_{L_{n}},$$
    and spherical balls $B_{n}$ with radius $\ell m_{d}/L_{n}$
    we have
\begin{align*}
    L_{n}^{-dq/2}  & \| P_{[L_{n}/\ell]}^{(1+\l,\l)}(\la b_{n},\cdot \ra) \|^{q}_{q}
    =L_{n}^{-dq/2} \| K_{n} \|_{p}^{p}
    \gtrsim
    \int_{\S^{d}}|L_{n}^{-d\ell /2}K_{n}(z)|^{p}d\mu_{L_{n}} (z)
    \\
    &
    \ge
    \int_{B_{n}}|L_{n}^{-d\ell /2}K_{n}(z)|^{p}d\mu_{L_{n}} (z)\ge \delta_{d}^{p}\mu_{L_{n}}(B_{n}).
\end{align*}
    and this together with (\ref{norma}) brings us the contradiction.

    Conversely, for any $z\in \S^{d}$
    and $Q\in \Pi_{L}$ we have
    $$|Q(z)|^{p}\le C_{d,\delta}L^{d+1}\int_{ \mathbb{B}(z,1/L) }|Q(u)|^{p}dm(u),$$
    where $\mathbb{B}(z,1/L)$ stands for the euclidean ball in $\R^{d+1}.$
    Using Corollary \ref{corl} we have
\begin{align*}
    \int_{\S^{d}}|Q(z)|^{p}d\mu_{L} (z) & \lesssim
    L^{d+1}\int_{\S^{d}}\int_{ \mathbb{B}(z,1/L) }|Q(u)|^{p}dm(u)d\mu_{L} (z)
    \\
    &
    \le
    L^{d+1}\int_{C_{1,L}}|Q(u)|^{p}\int_{\S^{d}}\chi_{\mathbb{B}(z,1/L)}(u)d\mu_{L} (z) dm(u)
    \\
    &
    \le
    L^{d+1}\int_{C_{1,L}}|Q(u)|^{p}\int_{\S^{d}}\chi_{B(u/|u|,1/L)}(z)d\mu_{L} (z) dm(u)
    \\
    &
    \le
    \frac{C}{\pi_{L}}L^{d+1}\int_{C_{1,L}}|Q(u)|^{p}dm(u)\sim \int_{\S^{d}}|Q(u)|^{p}d\s (u).
\end{align*}
\qed

\begin{corollary}                                                     \label{finite-union}
    Let $1\le p <\infty.$ The family $\mathcal{Z}\subset \S^{d}$ is a finite union of uniformly separated
    families
    if and only if there exists $C_{p}>0$ such that for all $L\ge 1$ and $Q\in \Pi_{L}$
\begin{equation}                                                            \label{ineq}
    \frac{1}{\pi_{L}}\sum_{j=1}^{m_{L}}|Q(z_{L j})|^{p}
    \le C_{p}\int_{\S^{d}}|Q(\omega)|^{p}d\sigma(\omega).
\end{equation}
\end{corollary}

\proof It is enough to take the family of measures
    $$\mu_{L}=\frac{1}{\pi_{L}}\sum_{j=1}^{m_{L}}\delta_{z_{Lj}},\;\;\; L\ge 0,$$
    and apply the previous result.
\qed

\begin{theorem}                                                                     \label{cont-unif-sep}
    Any $L^{p}-$MZ family $\mathcal{Z}$
    contains a uniformly separated family
    $\tilde{\mathcal{Z}}\subset \mathcal{Z}$ which is also an $L^{p}-$MZ family.
\end{theorem}

\proof
    First consider $1\le p<\infty.$
    Using Corollary \ref{finite-union} we can assume that $\mathcal{Z}$ is a finite union of $N$
    uniformly $\e-$separated families, that we call
    $\mathcal{Z}^{(j)},$ $j=1,\dots ,N.$
    Now, following \cite[p. 141]{Sei95} we can construct for $0<\delta<\e/4$ a uniformly separated
    family
    $\tilde{\mathcal{Z}}\subset \mathcal{Z}$ such that for all $L\ge 0$ and $j=1,\dots ,m_{L}$
    $$d(z_{Lj},\tilde{\mathcal{Z}}(L))<\delta/L.$$
    Let $\tilde{z}$ be the closest point in $\tilde{\mathcal{Z}}(L)$ to $z\in \mathcal{Z}(L).$
    Given $Q\in \Pi_{L}$
    there exists $z^{'}\in \R^{d+1}$ in the segment joining $z$ and $\tilde{z}$ such that
    $$|Q(z)-Q(\tilde{z})|\le |\nabla Q(z^{'})||z-\tilde{z}|\le \frac{\delta}{L}|\nabla Q(z^{'})|.$$
    Differentiating Poisson's formula
    $$Q(v)=\frac{1}{\s(\S^{d})}\int_{\partial B(z^{'},r)}\frac{r^{2}-|v-z^{'}|^{2}}{r|u-v|^{d+1}}Q(u)d\sigma_{r}(u),$$
    and evaluating in $z^{'}$ we obtain
    $$|\nabla Q(z^{'})|^{p}r^{d+p}\le C \| Q \|_{L^{p}(\partial B(z^{'},r))}^{p}$$
    where $C$ only depends on $p$ and $d.$
    Integrating with respect to $r$ in $[0,\e/ 2 L]$ we get
    $$|\nabla Q(z^{'})|^{p}\le C_{\e} L^{d+p+1}\int_{\mathbb{B}(z^{'},\e/ 2 L)}|Q(v)|^{p}dm(v).$$
    Observe that the balls $\mathbb{B}(z^{'},\e/ 2 L)$ are mutually
    disjoint therefore
\begin{align*}
    \| Q
    \|_{L^{p}(\S^{d})}^{p} & \sim \frac{1}{\pi_{L}}\sum_{z\in \mathcal{Z}(L)}|Q(z)|^{p}
    \lesssim
    \frac{1}{\pi_{L}}\sum_{j=1}^{N}\sum_{z\in \mathcal{Z}^{(j)}(L)}(
    |Q(z)-Q(\tilde{z})|^{p}+|Q(\tilde{z})|^{p})
    \\
    &
    \le
    \frac{1}{\pi_{L}}\sum_{j=1}^{N}
    \d^{p}L^{d+p}\int_{C_{\e/2,L}}|Q(v)|^{p}dm(v)+
    \frac{CN}{\pi_{L}}\sum_{z\in \tilde{\mathcal{Z}}(L)}|Q(z)|^{p}
    \\
    &
    \lesssim
    C_{\e,p,d,N}\delta^{p}\| Q\|_{L^{p}(\S^{d})}^{p}+\frac{CN}{\pi_{L}}\sum_{z\in
    \tilde{\mathcal{Z}}(L)}|Q(z)|^{p}.
\end{align*}
    We finish by taking $\delta$ small enough. The reverse inequality
    follows from Corollary \ref{finite-union}.

    For $p=\infty$ take $\e>0$ such that $2 C \e <1,$ where $C$ is the constant in the MZ inequality.
    Let $u,v \in \S^{d}$ be such that $d(u,v)<\e/L.$ Bernstein's inequality
    for trigonometric polynomials
    applied to the restriction of $Q$ to a great circle gives us
    $$|Q(u)-Q(v)|\le \e \| Q\|_{\infty},$$
    for $Q\in \Pi_{L}.$
    Now it is easy to construct a $\tilde{\Z}(L)\subset \Z(L)$ such that $d(u,v)>\e/L$
    for $u,v\in \tilde{\Z}(L)$ and any $z\in \Z(L)$
    belongs to a ball of center one point in $\tilde{\Z}(L)$ and radius $\e/L.$
    We denote $\tilde{\Z}(L)=\{ z_{Lk_{j}} \}_{j=1,\dots, N}$
    and for $Q\in \Pi_{L}$
\begin{align*}
    \| Q \|_{\infty} & \le C\sup_{z\in \Z(L)}|Q(z)|=C\max_{j=1,\dots, N}\sup_{z\in \Z(L),d(z,z_{Lk_{j}})<\e/L}|Q(z)|
    \\
    &
    \le C \e \| Q \|_{\infty}+C \max_{z\in \tilde{\Z}(L)} |Q(z)|.
\end{align*}
    So we obtain a $\e-$uniformly separated family $\tilde{\Z}$ such that for $Q\in \Pi_{L}$
    $$\| Q \|_{\infty}\le 2C \max_{z\in \tilde{\Z}(L)} |Q(z)|.$$
\qed

    Proposition \ref{prop} works only when $p=2.$ For other $p\in [1,\infty]$ we
    use a perturbative result.

\begin{definition}                                                                          \label{pertur}
    Given a family
    $\mathcal{Z}$ and $\delta>0$, we denote by
    $\mathcal{Z}_{\delta}$ the family
    $\mathcal{Z}_{\delta}(L)=\mathcal{Z}(L_{1+\delta}),$
    where $L_{1+\delta}=[(1+\delta)L].$
\end{definition}

\begin{lemma}                                                                       \label{rec-2}
    Let $p\in [1,\infty]$ and
    $\mathcal{Z}$ be a uniformly separated $L^{p}-$MZ family, then for $\delta>0$ and
    $q\in [1,\infty]$ the family
    $\mathcal{Z}_{\delta}$ is $L^{q}-$MZ.
\end{lemma}

\proof
    Using Riesz-Thorin theorem on interpolation of operators, see \cite[p.524]{DS},
    it is enough to show that $\mathcal{Z}_{\delta}$ is an $L^{q}-$MZ family for $q=1,\infty.$
    Fixed $z\in \S^{d}$ the evaluation operator $e_{z}(Q_{L})=Q_{L}(z)$ defined in $(\Pi_{L},\|\cdot \|_{p})$
    can be written as
    $$e_{z}(Q_{L})=\frac{1}{\pi_{L}}\sum_{j=1}^{m_{L}}Q_{L}(z_{Lj})a_{Lj}(z),$$
    where $a_{Lj}(z)\in \C$ are such that
    $\sum_{j=1}^{m_{L}}|a_{Lj}(z)|^{p^{'}}<C\pi_{L},$ where $1/p+1/p^{'}=1.$
    Let $p_{L_{\delta}}(t)$ be a polynomial in one variable of degree $L_{\delta}$
    such that $p_{L_{\delta}}(1)=1$ and
    $$\int_{-1}^{1}|p_{L_{\delta}}(t)|^{p}(1-t^{2})^{\lambda} dt=1.$$
    We have
    $$Q_{L}(z)=\frac{1}{\pi_{L_{1+\delta}}}
    \sum_{j=1}^{m_{L_{1+\delta}}}Q_{L}(z_{L_{1+\delta}j})p_{L_{\delta}}(z\cdot z_{L_{1+\delta}j})
    a_{L_{1+\delta}j}(z),$$
    so
\begin{align*}
    |Q_{L}(z)|\le & C \sup_{j}|Q_{L}(z_{L_{1+\delta}j})|
    \left( \frac{1}{\pi_{L_{1+\delta}}}
    \sum_{j=1}^{m_{L_{1+\delta}}} |p_{L_{\delta}}(z\cdot z_{L_{1+\delta}j})|^{p}\right)^{1/p}
    \\
    \le
    &
    C \left( \int_{-1}^{1}|p_{L_{\delta}}(t)|^{p}(1-t^{2})^{\lambda} dt\right)^{p} \sup_{j}|Q_{L}(z_{L_{1+\delta}j})|.
\end{align*}
    For $q=1$ we take
    $p_{L_{\delta}}(t)$ polynomial of degree $L_{\delta}$ in one variable
    such that $p_{L_{\delta}}(1)=1$ and
    $$\int_{-1}^{1}|p_{L_{\delta}}(t)|(1-t^{2})^{\lambda} dt=\pi_{L_{1+\delta}}^{-1}.$$
    and we get the result
\begin{align*}
    \int_{\S^{d}} |Q_{L}(z)|d\sigma(z)\le C
    \sum_{j=1}^{m_{L_{1+\delta}}}|Q_{L}(z_{L_{1+\delta}j})|
    \int_{\S^{d}} |p_{L_{\delta}}(z\cdot z_{L_{1+\delta}j})|d\sigma(z).
\end{align*}
\qed

    Finally we prove the corresponding result for interpolation. But first we want to estimate the norm of the
    evaluation operator.
    As in the proof of Theorem \ref{finite-union} we have for
    $Q\in \Pi_{L}$ and $u\in \S^{d}$
\begin{align*}
    |Q(u)|^{p} & \lesssim L^{d+1}\int_{B(u,1/L)}|Q(v)|^{p}d\sigma(v)
        \le
    L^{d+1}\int_{C_{1,L}}|Q(v)|^{p}dm(v)
    \\
    &
    \sim L^{d}\int_{\S^{d}}|Q(v)|^{p}d\sigma (v),
\end{align*}
    so $$\pi_{L}^{-1}\|Q\|_{\infty}^{p}\lesssim \| Q \|_{p}^{p}.$$

\begin{theorem}                                                                     \label{interp-unif}
    If $\mathcal{Z}$ is an interpolation family for $L^{p},$ then it is uniformly separated.
\end{theorem}

\proof
    Standard arguments based on the open mapping theorem for Banach spaces, see \cite{Sei95}, show that
    the interpolation can be done with polynomials
    $P_{L}$ such that
    $$|| P_{L}
    ||^{p}\lesssim \frac{1}{\pi_{L}}\sum_{j=0}^{m_{L}}|P_{L}(z_{Lj})|^{p}.$$
    Then, for a given $L_{0}\ge 0$ and $1\le j_{0}\le \pi_{L_{0}},$
    we can take polynomials $P_{L_{0}j_{0}}\in \Pi_{L_{0}}$ such that
    $P_{L_{0}j_{0}}(z_{Lj})=\delta_{LL_{0}}\delta_{jj_{0}}$
    and
    $\| P_{L_{0}j_{0}} \|^{p}_{p}\lesssim
    \pi_{L}^{-1}.$
    Then for
    $j\neq j_{0}$ restricting the polynomial to a great circle and using Bernstein's inequality for trigonometric polynomials
\begin{align*}
    1= & |P_{L_{0}j_{0}}(z_{L_{0}j_{0}})-P_{L_{0}j_{0}}(z_{L_{0}j})|\le
    \sup_{\gamma}|D_{T}P_{L_{0}j_{0}}| d(z_{L_{0}j_{0}},z_{L_{0}j})
    \\
    &
    \le L_{0} || P_{L_{0}j_{0}}
    ||_{\infty}d(z_{L_{0}j_{0}},z_{L_{0}j})\lesssim L_{0} \pi_{L}^{1/p} || P_{L_{0}j_{0}}
    ||_{p}d(z_{L_{0}j_{0}},z_{L_{0}j})
    \\
    &
    \lesssim
    L_{0} d(z_{L_{0}j_{0}},z_{L_{0}j}),
\end{align*}
    where $D_{T}$ stands for any unitary tangential derivative.
\qed

\begin{lemma}                                                                       \label{interp-ext}
    Let $p\in [1,\infty]$ and let
    $\mathcal{Z}$ be an $L^{p}-$interpolation family. For $\delta>0$ and
    $q\in [1,\infty]$
    $\mathcal{Z}_{-\delta}$ (as in Definition \ref{pertur}) is an $L^{q}-$interpolation family.
\end{lemma}

\proof
    As in the previous Lemma we will show that $\mathcal{Z}_{-\delta}$ is an $L^{q}-$interpolation family for $q=1,\infty.$
    The hypothesis implies that
    there exist polynomials $Q_{L_{1-\delta},j} \in \Pi_{L_{1-\delta}}$ such that
    $$Q_{L_{1-\delta},j}(z_{L_{1-\delta},k})=\delta_{jk},\;\;\;1\le j,k\le m_{L_{1-\delta}},$$
    with $$\| Q_{L_{1-\delta},j} \|_{p}^{p}\lesssim \pi_{L_{1-\delta}}^{-1}.$$
    Now take polynomials
    $p_{L_{\delta}}$ in one variable of degree $L_{\delta},$
    such that $p_{L_{\delta}}(1)=1$
    $$\int_{-1}^{1}|p_{L_{\delta}}(t)|^{p^{'}}(1-t^{2})^{\lambda} dt=\pi_{L_{1-\delta}}^{-1},
    \;\;\mbox{for}\;\;\frac{1}{p}+\frac{1}{p^{'}}=1.$$
    Given a triangular family $\{ c_{L_{1-\delta}j} \}_{L,j}$ such that
    $$\frac{1}{\pi_{L_{1-\delta}}}\sum_{j=1}^{m_{L_{1-\delta}}}
    |c_{L_{1-\delta}j}|<C,$$
    construct the polynomial
    $$Q_{L}(z)=\sum_{j=1}^{m_{L_{1-\delta}}} c_{L_{1-\delta}j}Q_{L_{1-\delta},j}(z)
    p_{L_{\delta}}(z\cdot z_{L_{1-\delta}j})\in \Pi_{L},$$
    which satisfies $Q_{L}(z_{L_{1-\delta}j})=c_{L_{1-\delta}j}$ and
\begin{align*}
    \int_{\S^{d}}|Q_{L}(z)|d\sigma (z) & \le
    \sum_{j=1}^{m_{L_{1-\delta}}} |c_{L_{1-\delta}j}|
    \| Q_{L_{1-\delta},j} \|_{p}
    \| p_{L_{\delta}}(\la \cdot, z_{L_{1-\delta}j}\ra) \|_{p^{'}}
    \\
    &
    \lesssim
    \frac{1}{\pi_{L_{1-\delta}}}\sum_{j=1}^{m_{L_{1-\delta}}}
    |c_{L_{1-\delta}j}|.
\end{align*}
    For $q=\infty$ we take polynomials $p_{L_{\delta}}$ as before, but with
    $$\int_{-1}^{1}|p_{L_{\delta}}(t)|(1-t^{2})^{\lambda} dt=\pi_{L_{1-\delta}}^{-1}.$$
    And defining $Q_{L}$ as before we obtain the interpolation property and
    $$|Q_{L}(z)|\le C \sup_{j}|c_{L_{1-\delta},j}|
    \sum_{j=1}^{m_{L_{1-\delta}}}|p_{L_{\delta}}(\la z,z_{L_{\delta}}\ra )|
    \le
    C \sup_{j}|c_{L_{1-\delta},j}|.$$
\qed

\section{There are no complete interpolation families in $L^{p}$ for $p\neq 2.$}                        \label{secnoCIS}
%==========================================================================================

    In this section we show that there are no $L^{p}-$complete interpolation
    families
    for $p\neq 2.$ We construct, using transference methods (see \cite[Theorem 1.1]{BC}), a projection
    in $L^{p}(\S^{d})$ that yields a
    bounded
    ball
    multiplier
    in $L^{p}(\R^{d}).$ Finally the celebrated result of C. Fefferman \cite{F}
    says that this can happen only for $p= 2.$

\proof (Theorem \ref{noCIS})
    We argue by contradiction. Let $\Z$ be an $L^{p}-$complete interpolation family.
    By Theorem \ref{interp-unif} we know that it is uniformly separated. Let $\e>0$ be the separation
    constant.
    Let $\ell^{p}_{L}$ be the vector space of $\{c_{j}\}\in \C^{L^{d}}$ with norm given by
    $\| \{c_{j}\} \|_{\ell^{p}_{L}}^{p}=\frac{1}{L^{p}}\sum_{j=1}^{L^{p}}|c_{j}|^{p}.$
    For $L\ge 0$ we consider the map $R_{L}:L^{p}(\S^{d})\longrightarrow \ell^{p}_{L}$ defined as
    $$L^{p}(\S^{d}) \ni f\longmapsto \{ \la f, L^{-d}K_{L}(\cdot ,z_{Lj})K_{2L}(\cdot ,z_{Lj})\ra
    \}_{j=1,\dots ,L^{d}}.$$
    We want to show that $R_{L}$ is bounded for $p=1,\infty,$ with constant independent of $L.$
    So let $f\in L^{1}(\S^{d}),$
\begin{align*}
    \frac{1}{L^{d}} \sum_{j=1}^{L^{d}} & |\la f,
    L^{-d}K_{L}(\cdot ,z_{Lj})K_{2L}(\cdot ,z_{Lj})\ra|
    \le \frac{1}{L^{d}}\sum_{j=1}^{L^{d}}
    \int_{\S^{d}}|f(\omega)|\left|
    \frac{K_{L}(\o , z_{Lj})}{L^{d}}K_{2L}(\o, z_{Lj})\right|d\omega
    \\
    &
    \le
    || f ||_{1}\frac{1}{L^{d}}\sup_{\omega\in\S^{d}}\sum_{j=1}^{L^{d}}\left|
    \frac{K_{L}(\o , z_{Lj})}{L^{d}}K_{2L}(\o, z_{Lj})\right|.
\end{align*}
    Let $\omega\in\S^{d}$ be fixed. Then
\begin{align*}
    \sum_{j=1}^{L^{d}} & \left|
    \frac{K_{L}(\o,z_{Lj})}{L^{d}}K_{2L}(\o,z_{Lj})\right|\sim
    \sum_{j=1}^{L^{d}}\left|
    P_{L}^{(1+\l,\l)}(\la z_{Lj}, \omega\ra )P_{2L}^{(1+\l,\l)}(\la z_{Lj}, \omega\ra )\right|
    \\
    &
    \le L^{d}+\sum_{j\in \mathcal{I}}\left|
    P_{L}^{(1+\l,\l)}(\la z_{Lj}, \omega\ra )P_{2L}^{(1+\l,\l)}(\la z_{Lj}, \omega\ra )\right|+L^{d-2}
\end{align*}
    where $\mathcal{I}$ are the indices $j$ such that $\frac{\e}{2(L+1)}\le d(\omega,z_{Lj})\le \pi-\frac{\e}{2(L+1)}.$
    Observe that there are only two points $z_{Lj}$ such that
    $j\not \in \mathcal{I}$ (one on each cap), and the value of the polynomial is bounded by the local maximum.
    In between we use Sz\"ego's estimate (\ref{estimate}) to get
\begin{align*}
    \sum_{j\in \mathcal{I}} & \left|
    P_{L}^{(1+\l,\l)}(\la z_{Lj}, \omega\ra )P_{2L}^{(1+\l,\l)}(\la z_{Lj}, \omega\ra )\right|\lesssim
    \frac{1}{L}\sum_{j\in \mathcal{I}}k^{2}(d(z_{Lj}, \o )).
\end{align*}

    Using rotation invariance we can suppose that
    $\omega=N.$ The function $k$ is decreasing in $(0,\pi/2)$ and a lot bigger around 0
    than around $\pi.$ Then to increase the sum we place the points $z_{Lj},$
    the closer the better,
    in "bands" around the north pole.
    Coarse estimates using the uniform separation
    yields $\# \mathcal{I}=O(L^{d}),$ a maximum of $O(L^{d-1})$ "bands"
    and $O(\frac{L}{\e}\sin \frac{\ell \e}{L})$ points in the $\ell-$th "band",
    if we start counting from $N.$
    So
    $$\frac{1}{L}\sum_{j\in \mathcal{I}}k^{2}(d(\omega,z_{Lj}))
    \lesssim \frac{1}{L}\sum_{\ell=1}^{L^{d-1}}\left( \frac{L}{\e}\right)^{d-1}\sin^{d-1}
    \frac{\ell \e}{L}k^{2}\left(\frac{\ell \e}{L}\right)
    \lesssim L^{d}\sum_{\ell=1}^{L^{d-1}}\frac{1}{\ell^{2}},$$
    and we get
    $$|| R_{L}f ||_{\ell^{1}_{L}}=\frac{1}{L^{d}} \sum_{j=1}^{L^{d}}|\la f,
    L^{-d}K_{L}( \cdot ,z_{Lj})K_{2L}( \cdot ,z_{Lj})\ra|\lesssim || f ||_{1},$$
    where the constat depends on $\e$ but is independent of $L.$
    To prove
    the $L^{\infty}$ case is a lot easier:
\begin{align*}
    L^{-d} & || K_{L}( \cdot ,z_{Lj})K_{2L}( \cdot ,z_{Lj}) ||_{1}
    \sim ||
    P_{L}^{(1+\l,\l)}(\la \cdot ,z_{Lj}\ra)
    P_{2L}^{(1+\l,\l)}(\la \cdot ,z_{Lj}\ra) ||_{1}
    \\
    &
    \le || P_{L}^{(1+\l,\l)}(\la \cdot ,z_{Lj}\ra) ||_{2}
    || P_{L}^{(1+\l,\l)}(\la \cdot ,z_{Lj}\ra) ||_{2}=|\S^{d}|.
\end{align*}
    Now let
    $E_{L}$ be the map from $\ell^{p}_{L}$ to
    $\Pi_{L},$ sending
    $v=\{ v_{j} \} \in \ell_{L}^{p}$ to
    $P_{L}\in
    \Pi_{L}$ such that $P_{L}(z_{Lj})=v_{j}.$ By hypothesis
    $|| E_{L}(v) ||_{p}\asymp\|v \|_{\ell^{p}_{L}},$
    so $E_{L}\circ R_{L}$ is bounded from
    $L^{p}(\S^{2})$ to $L^{p}(\S^{2})$ for
    $p=1,\infty$ and
    by Riesz-Thorin theorem on interpolation of operators, see \cite[p.524]{DS}, we get that it is bounded
    for all $1\le p\le \infty.$
    Denoting $\mathcal{P}_{L}=E_{L}\circ R_{L}$ we get
    ${\mathcal{P}_{L}}_{|\Pi_{L}}=I_{\Pi_{L}}.$

    Following \cite[Theorem 1]{R} we define
    $$\mathfrak{P}_{L}f=\int_{SO(d+1)}\nu^{-1}\mathcal{P}_{L}\nu f
    d\nu,$$
    that turns out to be a projection from
    $L^{p}(\S^{d})$ to $\Pi_{L},$
    commuting with rotations and such that $|| \mathfrak{P}_{L} ||\le
    || \mathcal{P}_{L} ||.$

    According to
    Theorem \ref{co_we}
    we have
    $\mathfrak{P}_{L}Y=m_{\ell}Y,$
    for $Y\in \mathcal{H}_{\ell}$ and
    for $m_{\ell}\in \C.$
    The properties of $\mathfrak{P}_{L}$ impose that
    $m_{\ell}=1$ for $\ell\le L$ and zero otherwise. So
    $\mathfrak{P}_{L}f$ is just the sum of the orthogonal projections of $f$ over
    $\mathcal{H}_{\ell}$ (denoted by $P_{\mathcal{H}_{\ell}}f$) for $\ell=0,\dots ,L.$

    Now we can put
    $$\mathfrak{P}_{L}f=\sum_{j=0}^{\infty}m_{L}(\ell)P_{\mathcal{H}_{\ell}}f,$$
    with $m_{L}(\ell)=m(\frac{\ell}{L})$ and
    $m(|x|)=\chi_{\mathbb{B}}(x).$ The sequence $\{m_{L}(\ell)\}_{\ell\ge 0}$
    defining a multiplier in $L^{p}(\S^{2})$ with
    $$\sup_{L\ge 0} ||\mathfrak{P}_{L} ||_{p}<\infty$$
    where $m_{0}(\ell)=\delta_{0\ell}.$

    Now using the transference result in \cite[Theorem 1.1]{BC}
    we see that the multiplier in  $L^{p}(\R^{d})$ given by
    $$f \longmapsto \mathcal{F}^{-1}(\chi_{\mathbb{B}}\mathcal{F}f),$$
    is bounded.
    Finally C. Fefferman's result \cite{F} says that this is only possible for
    $p=2.$
\qed

\section{Proofs}                                      \label{SecProofs}
%--------------------------------------------------------------------------------------

    We need some notation and two technical Lemmas before proving Theorem \ref{teo}.

    Given $L\ge 0$ and $\alpha>0$
    let $A_{L},$ $A_{L}^{+}$ and $A_{L}^{-}$ be
    the geodesic balls centered at the north pole
    with respective radius $\alpha/(L+1),$ $(\alpha+\e)/(L+1)$ and $(\alpha-\e)/(L+1),$
    where $\e$ will denote the separation constant.

    Denote the eigenvalues of
    the concentration operator
    $\K_{A_{L}}$
    as
    $$1>\lambda_{1}^{L}\ge \dots \ge \lambda_{\pi_{L}}^{L}>0.$$

\begin{lemma}                                                       \label{lema1}
    Let  $\mathcal{Z}$ be a $\e-$uniformly separated
    $L^{2}-$MZ family and let
    $$N_{L}=\# (\mathcal{Z}(L)\cap A_{L}^{+}).$$
    There exists a constant $0<\gamma<1$ independent of
    $\alpha$ and $L$ such that
    $$\lambda_{N_{L}+1}^{L}\le \gamma.$$
\end{lemma}

\remark
    In the conditions of the Lemma \ref{lema1}
    $$
    \#\{\l_{j}^{L}>\gamma \} \le N_{L}= \# (\mathcal{Z}(L)\cap A_{L}^{+})
    \le
    \# (\mathcal{Z}(L)\cap
    A_{L})+C (1+o(\a^{d})),\;\;\;\alpha\to\infty,$$
    where the constant $C$ depends on $d$ and $\e.$
    This follows from the estimates
    $L^{d}\sigma(A_{L}^{+}\setminus A_{L})=1+o(\a^{d})$
    if $\a\to\infty$
    and
    $$\# (\mathcal{Z}(L)\cap (A_{L}^{+}\setminus A_{L}))
    \frac{\e^{d}}{L^{d}}\lesssim \sigma(A_{L}^{+}\setminus A_{L}).$$

\begin{lemma}                                                                   \label{lema2}
    Let $\mathcal{Z}$ be an $L^{2}-$interpolation family and let
    $$n_{L}=\# (\mathcal{Z}(L)\cap A_{L}^{-}).$$
    There exists a constant $0<\delta<1$ independent of
    $\alpha$ and $L$
    such that
    $$\lambda_{n_{L}-1}^{L}\ge \delta.$$
\end{lemma}

\remark In the conditions of the Lemma \ref{lema2}
    we have, as before,
    $$
    \#(\mathcal{Z}(L)\cap
    A_{L})-C(1+o(\a^{d}))
    \le n_{L}=\# \mathcal{Z}(L)\cap
    A_{L}^{-}
  \le
    \#\{\l_{j}^{L}\ge \delta\}+1.$$

\proof[Theorem \ref{teo}]
    Using Theorems \ref{cont-unif-sep} and \ref{interp-unif} we can
    suppose that $\Z$ is
    a uniformly separated family.
    Now given $\eta>0$ and taking either $\Z_{\eta}$ or $\Z_{-\eta}$
    we have
    by Lemmas \ref{interp-ext} and \ref{rec-2} that our family is respectively $L^{2}-$MZ or interpolating.
    Now we relabel the family as before and
    defining the measures $d\mu_{L}=\sum_{j=1}^{\pi_{L}}\delta_{\l_{j}^{L}}$ we
    have
    $$\tr (\mathcal{K}_{A_{L}})=\int_{0}^{1}xd\mu_{L}(x),\;\;\;\;\mbox{and}
    \;\;\;\;\tr (\mathcal{K}_{A_{L}}^{2})=\int_{0}^{1}x^{2}  d\mu_{L}(x).$$
    Let $\mathcal{Z}$ be an $L^{2}-$MZ and
    let $\gamma$ be given by Lemma \ref{lema1}. We get
\begin{align*}
    \# \{ \l_{j}^{L}>\gamma \} & =
    \int_{\gamma}^{1}d\mu_{L}(x)
    \ge
    \int_{0}^{1}x d\mu_{L}(x)
    -\frac{1}{1-\gamma}\int_{0}^{1}x(1-x)d\mu_{L}(x)
    \\
    &
    =
    \tr (\K_{A_{L}})-\frac{1}{1-\gamma}(\tr (\K_{A_{L}})-\tr (\K_{A_{L}}^{2})).
\end{align*}
    The remark following Lemma \ref{lema1} and Proposition \ref{prop}
    yield
\begin{align*}
    \frac{\# ( \mathcal{Z}(L)\cap
    A_{L})+C(1+o(\a^{d}))}{\alpha^{d}}
    \ge
    \frac{\pi_{L}\s (A_{L})}{\a^{d}\s (\S^{d})}
    %\frac{2|\S^{d-1}|}{|\S^{d}|d!d}
    &
    -\frac{O(\a^{d-1}\log \a)}{\alpha^{d}(1-\gamma)},
\end{align*}
    and taking limits
    we get,
    for any $\eta>0,$
    $$D^{-}(\Z_{\eta})\ge
    %\frac{2}{\sqrt{2\pi}(d!!)^{2}},$$
    \frac{2}{d!d\sqrt{\pi}}\frac{\Gamma(\frac{d+1}{2})}{\Gamma(\frac{d}{2})},$$
    what implies the result.

    Assume now that $\mathcal{Z}$ is an $L^{2}-$interpolation family and
    let $\delta>0$ be the value provided by Lemma \ref{lema2}. Using the estimate of
    Proposition \ref{prop} we get
\begin{align*}
    \# & \{ \l_{j}^{L}\ge \delta \} \le
    \frac{-1}{\delta}\tr(\K^{2}_{A_{L}})+\frac{1+\delta}{\delta}\tr (\K_{A_{L}})
    \\
    &
    =
    \tr (\K_{A_{L}})+\frac{1}{\delta}(\tr (\K_{A_{L}})-\tr (\K^{2}_{A_{L}}))
    =
    \frac{\pi_{L}\s (A_{L})}{\s (\S^{d})}+\frac{1}{\delta}O(\a^{d-1}\log \a).
\end{align*}
    Using as before the remark following Lemma \ref{lema2}
    and taking limits
    we get
    for any $\eta>0$
    $$D^{+}(\Z_{\eta})\le
     \frac{2}{d!d\sqrt{\pi}}\frac{\Gamma(\frac{d+1}{2})}{\Gamma(\frac{d}{2})},$$
    %\frac{2\s(\S^{d-1})}{\s(\S^{d})d!d},$$
    what finishes the proof.
\qed

%demostracio lema 1------------------------------------------------------------------

    In the proof of the Lemmas \ref{lema1} and \ref{lema2} we follow
    \cite{L}. For the definition of the Gegenbauer polynomials and related notions see \cite{Mu}.

    Given
    $\delta>0$
    consider the functions
\begin{equation}                                                                            \label{funcioh}
    h(\omega)=\left(\frac{L}{\d} \right)^{d}\chi_{B(N,\frac{\delta}{2(L+1)})}(\omega),\;\;\;\omega \in \S^{d}.
\end{equation}
    The polynomial $Y_{\ell}^{1}$ (a multiple of the Legendre harmonic)
    is just the Gegenbauer polynomial $C_{\ell}^{\frac{d-1}{2}}$ normalized in the $L^{2}-$norm.
    Applying Funk-Hecke theorem to $h$ we get
\begin{align*}
    \hat{h}(\ell,1) & =
    \left(\frac{L}{\d} \right)^{d}\int_{\S^{d}}\chi_{(\cos \frac{\delta}{2(L+1)},1)}(\la \omega, N\ra )
    Y_{\ell}^{1}(\omega)d\sigma(\omega)
    \\
    &
    =
    \frac{L^{d}\s(\S^{d-1})}{\d^{d}C_{\ell}^{\frac{d-1}{2}}(1)\|  C_{\ell}^{\frac{d-1}{2}}(\la N,\cdot \ra ) \|_{2}}
    \int_{0}^{\frac{\delta}{2(L+1)}}
     C_{\ell}^{\frac{d-1}{2}}(\cos \theta)\sin^{d-1}\theta \,d\theta .
    \end{align*}
    Given $f\in L^{2}(\S^{d})$ $0\le
    \ell,$ $1\le m\le h_{\ell}$ and applying Funk-Hecke as before,
    we deduce that
\begin{align*}
    (f \ast h)^{\hat{}} & (\ell,m)= \int_{\S^{d}}(f \ast h)(\omega)\overline{Y_{\ell}^{m}(\omega)}d\sigma(\omega)
     \\
     &
     =\frac{L^{d}}{\d^{d}}\int_{\S^{d}}f(u)\left( \int_{\S^{d}}\chi_{(\cos \frac{\d}{2(L+1)},1)}(\la u ,\o\ra )
    \overline{Y_{\ell}^{m}(\omega)}d\s (\o)\right)
     d\s (u)
    \\
    &
    =
    \|  C_{\ell}^{\frac{d-1}{2}}(\la N,\cdot \ra ) \|_{2}\hat{h}(\ell,1)
    \hat{f}(\ell,m)
\end{align*}
    thus
    $$|(f \ast h )^{\hat{}}(\ell,m)|=\frac{C_{L,\delta}\s(\S^{d-1})}{C_{\ell}^{\frac{d-1}{2}}(1)\s(\S^{d})}
    |\hat{f}(\ell,m)|\left|
    \int_{0}^{\frac{\delta}{2(L+1)}}C_{\ell}^{\frac{d-1}{2}}(\cos \theta)\sin^{d-1} \theta \,d\theta\right|.$$
    Now we want to show that
    for $0\le \ell\le L$ and $\delta$ sufficiently small
\begin{equation}                                                                        \label{desi}
    \left|
    \int_{0}^{\frac{\delta}{2(L+1)}}C_{\ell}^{\frac{d-1}{2}}(\cos \theta)\sin^{d-1} \theta
    d\theta\right|\gtrsim C_{\ell}^{\frac{d-1}{2}}(1)\left(\frac{\delta}{L}\right)^{d},
\end{equation}
    and in particular
    for all $Q\in \Pi_{L}$
    $$|(Q \ast h )^{\hat{}}(\ell,m)|\gtrsim |\hat{Q}(\ell,m)|,\;\;\;0\le \ell\le L\;\;1\le m\le h_{\ell}.$$

    To prove (\ref{desi}) let $x_{\ell}$ be the largest zero
    in $[-1,1]$ of $C_{\ell}^{\frac{d-1}{2}}.$ It is known that
    $x_{\ell}\sim \cos C/L,$ for some constant $C>0,$
    so for $\delta$ sufficiently small independent of L, the polynomial
    $C_{\ell}^{\frac{d-1}{2}}$ has no zeros in
    the spherical cap centered in $N$ with radius $\delta/2(L+1),$ \cite{S}.
    The integral in (\ref{desi}) can be written as
    $$\int_{B(N,\frac{\delta}{2(L+1)})}C_{\ell}^{\frac{d-1}{2}}(\la \o ,
    N\ra )d\sigma(\omega),$$
    and for $d(\omega,N)<\delta/2(L+1)$
    $$C_{\ell}^{\frac{d-1}{2}}(\la \o ,
    N\ra )\ge C_{\ell}^{\frac{d-1}{2}}(1)\left(
    1-\frac{2(L+1)d(\omega,N)}{\delta}\right),$$
    or equivalently
    $$C_{\ell}^{\frac{d-1}{2}}(x)\ge
    C_{\ell}^{\frac{d-1}{2}}(1)\left( 1-\frac{2(L+1)\arccos
    x}{\delta}\right),$$
    if $\cos \frac{\delta}{2(L+1)}\le x\le 1.$ This can be deduced using the concavity of the polynomial
    and the convexity of the function in the right hand side of the last expression.
    So
\begin{align*}
    \int_{B(N,\frac{\delta}{2(L+1)})} & C_{\ell}^{\frac{d-1}{2}}(\omega\cdot
    N)d\sigma(\omega)
    \ge \int_{B(N,\frac{\delta}{2(L+1)})}C_{\ell}^{\frac{d-1}{2}}(1)\left(
    1-\frac{2(L+1)d(\omega,N)}{\delta}\right)d\sigma(\omega)
    \\
    &
    \sim
    \int_{0}^{\frac{\delta}{2(L+1)}}C_{\ell}^{\frac{d-1}{2}}(1)\sin^{d-1}\eta
    \left(
    1-\frac{2(L+1)\eta}{\delta}\right)d\eta
    \\
    &
    =
    C_{\ell}^{\frac{d-1}{2}}(1)\int_{0}^{1}\frac{\delta}{2(L+1)}\sin^{d-1}\left(\frac{\delta}{2(L+1)}(1-\eta)\right)
    d\eta
    \\
    &
    \gtrsim
    C_{\ell}^{\frac{d-1}{2}}(1)\left(
    \frac{\delta}{2(L+1)}\right)^{d},
\end{align*}
    and (\ref{desi}) follows.
\qed

%%demostracio_lema1---------------------------------------------------------------------------------------

\proof (Lemma (\ref{lema1}))
    Let $Q \in \Pi_{L},$ let $0<\delta<\e,$ where $\e>0$ is the separation constant of $\mathcal{Z}$
    and let $h$ be as in (\ref{funcioh}).
    Defining
    $g=Q\ast h\in \Pi_{L}$ we have
    $$||Q||^{2}
    =\sum_{\ell=0}^{L}\sum_{k=1}^{h_{\ell}}|\hat{Q}(\ell,k)|^{2}\lesssim
    \sum_{\ell=0}^{L}\sum_{k=1}^{h_{\ell}}|(Q\ast h)^{\hat{}}(\ell,k)|^{2}=|| g ||^{2}
    \lesssim
    \frac{1}{\pi_{L}}\sum_{k=1}^{m_{L}}|g(z_{Lj})|^{2}.$$
    Applying Schwarz's inequality
    $$
    |g(z_{Lj})|^{2} =\left| \int_{\nu \in SO(d+1)}Q(\nu N)h(\nu^{-1}z_{Lj})d\nu\right|^{2}
    \le
    \frac{|| h ||^{2}_{L^{2}(\S^{d})}}{\s(\S^{d})} \int_{d(\nu N,z_{Lj})<\frac{\epsilon}{2(L+1)}}|Q(\nu N)|^{2}d\nu.$$
    Now suppose that
    $$g(z_{Lj})=0, \;\;\mbox{for any}\;\; z_{Lj}\in A_{L}^{+},$$
    and denote by $\mathcal{I}$ the set of indices of those points $z_{Lj}$ where $g$
    vanishes.
    Then
\begin{align*}\label{uno}
    || Q ||^{2} & \lesssim \frac{1}{\pi_{L}}\sum_{j \not \in \mathcal{I}}|g(z_{Lj})|^{2}
    \le
    \frac{|| h ||^{2}_{L^{2}(\S^{d})}}{\pi_{L}}\sum_{j \not \in \mathcal{I}}
    \int_{d(\nu N,z_{Lj})<\frac{\epsilon}{2(L+1)}}|Q(\nu N)|^{2}d\nu
    \\
    &
    \le
    C_{\delta} \int_{\S^{d}\setminus A_{L}}|Q(\omega)|^{2}d\sigma(\omega),
\end{align*}
    where we have used the separation in the last inequality.

    Now we consider an orthonormal basis of eigenvectors $G_{j}^{L},$ corresponding to the eigenvalues
    $\lambda_{j}^{L}$ and let $c_{j}^{L}$ in
    $$Q(z)=\sum_{j=1}^{N_{n}+1}c_{j}^{L}G_{j}^{L}\in \Pi_{L},$$ be such that
    $g(z_{Lj})=(Q\ast h)(z_{Lj})=0$ for $z_{Lj}\in A_{L}^{+}.$

    Now
\begin{align*}
    \lambda_{N_{L}+1}^{L} & \sum_{j=0}^{N_{L}+1}|c_{j}^{L}|^{2}\le \sum_{j=0}^{N_{L}+1}\lambda_{j}^{L}|c_{j}^{L}|^{2}
    =|| \chi_{A_{L}}Q||^{2}=||Q ||^{2}-|| \chi_{\S^{2}\setminus A_{L}}Q||^{2}
    \\
    &
    \le \left( 1-\frac{1}{C_{\delta}}\right)\sum_{j=0}^{N_{L}+1}|c_{j}^{L}|^{2},
\end{align*}
    and we get the result.
\qed

%%demostracio lema2--------------------------------------------------------------------------------------

\proof (Lemma \ref{lema2})
    Let $\widetilde{\Pi}_{L}$ be the subspace of those polynomials in $\Pi_{L}$
    vanishing in $\mathcal{Z}(L).$
    Let $Q_{j}\in\Pi_{L}\ominus \widetilde{\Pi}_{L}$ be such that
    $$Q_{j}(z_{Lj^{'}})=\delta_{j j^{'}},$$
    and let
    $h$ be as in (\ref{funcioh})
    with
    $0<\delta<\e$ where
    $\e>0$ is the separation constant of $\mathcal{Z}.$

    Let $\tilde{Q_{j}}\in
    \Pi_{L}$ be such that
    $Q_{j}(\omega)=(\tilde{Q}_{j}\ast h)(\omega),$ and for
    $$Q\in \mbox{span}\{\tilde{Q}_{j}:z_{Lj}\in A_{L}^{-}
    \}$$
    we take $g=Q\ast h.$

    It is clear that $g\in \Pi_{L}\ominus \widetilde{\Pi}_{L}$ and vanishes in those points such that
     $z_{Lj}\not\in A_{L}^{-}.$
     Now following the same steps of
    Lemma \ref{lema2} and using that $g\in \Pi_{L}\ominus \widetilde{\Pi}_{L}$
    we get
    $$|| Q ||^{2}\le C\int_{A_{L}}|Q(\omega)|^{2}d\sigma(\omega).$$
    Applying Weyl-Courant's Lemma, \cite{DS},
    $$\l_{k-1}^{L}\ge \inf_{Q\in \Pi_{L},Q\in E}\frac{|| \chi_{A_{L}}Q ||^{2}}{|| Q
    ||^{2}}, \;\; \mbox{if}\;\;\dim E=k.$$
    Taking $E=\mbox{span}\{\tilde{Q}_{j}:z_{Lj}\in A_{L}^{-}
    \},$ that has dimension $n_{L},$
    we get the result.
\qed

%%---------------------------------------------------------------------------------------------------

%%-----------------------------------------------------------------------------------------------------------

\end{document}